\documentclass[a4paper]{amsart}
\usepackage[all]{xy}
\usepackage{amssymb}
\usepackage{amsmath}
\usepackage{bm}
\usepackage[colorlinks,linkcolor=blue]{hyperref}
\usepackage{color,xcolor,dsfont}
\usepackage{enumerate}
\usepackage{epsfig}
\usepackage{extarrows}
\usepackage{framed}
\usepackage{graphicx}
\usepackage{marvosym}
\usepackage{mathrsfs}
\usepackage{setspace} 
\usepackage{tikz}
\usepackage{tikz-cd}

\newcommand{\scha}{\mathcal{A}}

\newcommand{\schc}{\mathcal{C}}
\newcommand{\schd}{\mathcal{D}}
\newcommand{\schf}{\mathcal{F}}
\newcommand{\schk}{\mathcal{K}}
\newcommand{\schp}{\mathcal{P}}

\newcommand{\scht}{\mathcal{T}}
\newcommand{\she}{\mathscr{E}}
\newcommand{\shf}{\mathscr{F}}
\newcommand{\shg}{\mathscr{G}}

\newcommand{\shi}{\mathscr{I}}

\newcommand{\sho}{\mathscr{O}}

\newcommand{\shu}{\mathscr{U}}

\newcommand{\deriver}{\textup{\textsf{R}}}
\newcommand{\derivel}{\textup{\textsf{L}}}

\DeclareMathOperator{\id}{id}
\DeclareMathOperator{\GL}{GL}
\DeclareMathOperator{\NS}{NS}
\DeclareMathOperator{\Pic}{Pic}
\DeclareMathOperator{\Bir}{Bir}
\DeclareMathOperator{\Aut}{Aut}

\DeclareMathOperator{\Hom}{Hom}

\DeclareMathOperator{\coho}{coho}

\DeclareMathOperator{\rank}{rank}
\DeclareMathOperator{\Stab}{Stab}

\newcommand{\stable}{\textsf{\textup{st}}}

\newcommand{\bz}{\mathbb{Z}}

\newcommand{\br}{\mathbb{R}}
\newcommand{\bc}{\mathbb{C}}

\newcommand{\bp}{\mathbb{P}}

\newcommand{\ch}{\textup{\textsf{ch}}}
\newcommand{\td}{\textup{\textsf{td}}}

\newcommand{\Gro}{\textup{\textsf{K}}}
\newcommand{\Num}{\textup{\textsf{N}}}

\newcommand{\GLp}{\operatorname{GL^+}}
\newcommand{\grp}{{\tilde{\GLp}}(2,\br)}

\newcommand{\defi}[1]{{\textbf{\emph{#1}}}}
\newcommand{\cate}[1]{{\textbf{\textup{#1}}}}
\newcommand{\funct}[1]{{\textup{\textbf{\textsf{#1}}}}}

\theoremstyle{plain}
\newtheorem{theorem}{Theorem}[section]

\newtheorem{proposition}[theorem]{Proposition}

\theoremstyle{definition}
\newtheorem{definition}[theorem]{Definition}
\newtheorem{remark}[theorem]{Remark}
\newtheorem{example}[theorem]{Example}

\author{Ziqi Liu}
\title[derived-natural automorphisms]{derived-natural automorphisms on Hilbert schemes of points on generic K3 surfaces}
\address{Dipartimento di Matematica “F. Enriques”, Università degli Studi di Milano, Via Cesare Saldini 50, 20133 Milano, Italy.}
\email{ziqi.liu@unimi.it}



\thanks{The author is a member of the INdAM group GNSAGA (2025).}

\begin{document}

\maketitle
\thispagestyle{plain}

\begin{abstract}
The article revisits birational and biregular automorphisms of the Hilbert scheme of points on a K3 surface from the perspective of derived categories. Under the assumption that the K3 surface is generic, the birational and biregular involutions induced by autoequivalences on the derived category of the underlying K3 surface are characterized. 
\end{abstract}

\section{Introduction}
Hyperkähler manifolds are higher dimensional analogues of K3 surfaces. The space $H^2(X,\bz)$ of a hyperkähler manifold $X$ carries a canonical polarized weight-two Hodge structure \cite{Bea83} which, according to the global Torelli theorem, largely controls its isomorphism class. In particular, one has two representations
$$\rho_A\colon \Aut(X)\rightarrow O(H^2(X,\bz),q_X)\quad\textup{and}\quad \rho_B\colon \Bir(X)\rightarrow O(H^2(X,\bz),q_X)$$
where $q_X$ is the polarization for the canonical Hodge structure on $H^2(X,\bz)$.

Only limited examples of hyperkähler manifolds are known up to deformation, which are all obtained from moduli spaces of semistable sheaves on symplectic surfaces. Among these known deformation prototypes, the Hilbert scheme $S^{[n]}$ of $n$ points on a K3 surface $S$ has the simplest structure and is thus the best understood. In particular, the kernels of $\rho_A$ and $\rho_B$ for it are both trivial according to \cite{Bea83,Hu99}. As a result, studying the groups $\Aut(S^{[n]})$ and $\Bir(S^{[n]})$ reduces to studying certain Hodge isometries on its second integral cohomology space.

Thanks to the tools provided by \cite{BM14-MMP} and many other prior works, these Hodge isometries can be worked out at least theoretically. However, it is still a difficult task to describe a biregular or birational automorphism on $S^{[n]}$, even when its Hodge isometry is understood. In order to give descriptions, a notion of \emph{natural} automorphisms on $S^{[n]}$ is proposed in \cite{Boi12}: an automorphism on $S^{[n]}$ is called natural if it is induced by an automorphism on $S$. The sufficient and necessary condition for an automorphism to be natural has been determined in \cite{BS12}.

On the other hand, that is far from the end of the story as natural automorphisms are rather special: while the group $\Aut(S)=\Bir(S)$ is trivial when $S$ is generic and with degree larger than $2$ (see e.g.\ \cite[Corollary 15.2.12]{HuK3}), there can be non-trivial elements in $\Aut(S^{[n]})$ or $\Bir(S^{[n]})$ by \cite{BC22,BCWS16,Ca19}. Therefore, to give a characterization for these automorphisms, we need to generalize the notion of a natural automorphism. In this article, we consider the automorphisms on $S^{[n]}$ induced by autoequivalences on the bounded derived category $\cate{D}^b(S)$ of $S$ instead of by automorphisms on $S$. Such an automorphism is called \emph{derived-natural}.

According to \cite[Proposition 3.5]{Ou18}, an automorphism $f$ on $S^{[n]}$ is derived-natural only if the induced Hodge isometry $f_*$ on $H^2(S^{[n]},\bz)$ acts on the discriminant group $A_{H^2(S^{[n]},\bz)}$ of $H^2(S^{[n]},\bz)$ as the identity. This immediately raises the question: is this necessary condition also a sufficient one? The analogous question for natural automorphisms is answered in \cite{BS12}. This article provides a partial answer.

\begin{theorem}\label{1-main-result}
Consider a generic polarized K3 surface $(S,H)$ and a biregular or birational automorphism $f$ on $S^{[n]}$. Then $f$ is derived-natural if and only if the Hodge isometry $f_*$ induces the identity map on the discriminant group $A_{H^2(S^{[n]},\bz)}$.
\end{theorem}

The possible derived-natural automorphisms on $S^{[n]}$ are characterized in \cite{BC22,Ca19} whose Hodge isometries have been described. The theorem is proved by showing that these isometries can be lifted to the corrected autoequivalences using \cite{Hu16}.

\begin{remark}
The non-derived-natural biregular automorphisms are studied in \cite{FMOS22} and some of them can be realized by anti-autoequivalences \cite{ASF15,FMP25,Mar01}. A similar but different construction has been discussed in \cite[Section 3.1]{Pl07}.
\end{remark}

Besides the existence, there are other natural problems to consider such as the uniqueness. In this article, one can see

\begin{proposition}\label{1-non-main-result}
Under the setups of Theorem \ref{1-main-result}, any two autoequivalences inducing the same derived-natural automorphism are isomorphic to each other. 
\end{proposition}

At the end of this article, we remark that the notion of derived-natural automorphisms on $S^{[n]}$ can be generalized to other hyperkähler manifolds.

\subsection{Conventions}
Everything is over the field $\bc$ of complex numbers, the real and imaginary parts of $c\in\bc$ are denoted by $\Re c$ and $\Im c$. A generic polarized K3 surface $(S,H)$ means a K3 surface $S$ with $\Pic(S)$ generated by an ample class $H$.

\section{Stability Conditions}
\subsection{Basic definitions}
We recall the notion of a stability condition in \cite{Bri07}.

\begin{definition}
A \defi{stability function} $Z\colon \Gro(\scha)\rightarrow\bc$ on an abelian category $\scha$ is a group homomorphism from its Grothendieck group $\Gro(\scha)$ to the additive group $\bc$ such that $\Im Z(E) \geq 0$, and in the case that $\Im Z(E) = 0$, we have $\Re Z(E)< 0$.
\end{definition}

The \emph{slope function} $\mu_Z$ associated with a given stability function $Z$ is defined by $\mu_Z(E)= -\Re Z(E)/\Im Z(E)$ for an object $E\in\scha$ with $\Im Z(E)>0$, and $\mu_Z(E)=+\infty$ when $\Im Z(E)=0$. A non-zero object $E$ in $\scha$ is said to be \emph{semistable} (resp.\ \emph{stable}) with respect to $Z$ if for every proper subobject $A\subset E$ one has $\mu_Z(A)\leq \mu_Z(E/A)$ (resp.\ $\mu_Z(A)<\mu_Z(E/A)$). The \emph{phase} $\phi_Z$ of an object $E\in\scha$ associated with $Z$ is defined by $\phi_Z(E)=(1/\pi)\arg Z(E)\in(0,1)$  when $\Im Z(E)>0$, and $\phi_Z(E)=1$ otherwise. The subscripts for $\mu_Z$ and $\phi_Z$ associated with a stability function $Z$ will usually be omitted if the stability function is clear from the context.

\begin{definition}
Let $Z$ be a stability function on an abelian category $\scha$, then a non-zero object $E$ of $\scha$ is said to admit a \defi{Harder--Narasimhan (HN) filtration} if there exists a filtration $0=E_0 \hookrightarrow E_1 \hookrightarrow \dots E_{n-1} \hookrightarrow E_n=E$ whose factors $A_i:=E_i/E_{i-1} \neq 0$ are semistable with respect to $Z$ and $\mu(A_1) > \dots > \mu(A_n)$. 
\end{definition}

\begin{remark}
Once a Harder--Narasimhan filtration of an object exists, it is unique up to isomorphisms. Objects admitting the same factors are called \emph{$S$-equivalent}.
\end{remark}

A triangulated category $\schc$ is called of \emph{finite type} (over $\bc$) if for every pair of objects $E$ and $F$ the vector space $\oplus_i\Hom(E,F[i])$ is finite dimensional. In this situation one can define the Euler form on $\Gro(\schc)$ via the formula
$$\chi(E,F)=\sum(-1)^i\dim\Hom_{\schc}(E,F[i])$$
and the \emph{numerical
Grothendieck group} $\Num(\schc) := \Gro(\schc)/\Gro(\schc)^{\perp}$.

\begin{definition}
Let $\schc$ be a finite type triangulated category, then a \defi{numerical stability condition} $\sigma$ on $\schc$ is a pair $(\scha, Z)$, where $\scha$ is the heart of a bounded $t$-structure on $\schc$ and $Z\colon\Gro(\schc)\rightarrow \bc$ is a stability function on $\scha$ such that $\Gro(\schc)^{\perp}\subset\ker(Z)$ and every non-zero object $E$ in $\scha$ admits a Harder--Narasimhan filtration with respect to $Z$. A non-zero object $E$ of $\schc$ is said to be \defi{$\sigma$-semistable} (resp.\ \defi{$\sigma$-stable}) if a shift of $E$ is semistable (resp.\ stable) with respect to $Z$. 
\end{definition}

The slope can only be defined for $\sigma$-semistable object in hearts, so it is often convenient to order $\sigma$-semistable objects in $\schc$ by their phase which is $\phi_{\sigma}(E)=\phi_Z(E[n])-n$ for a $\sigma$-semistable object $E$ with $E[n]$ being in the heart $\scha$.

One can define for each $\phi\in\br$ an additive subcategory $\schp(\phi)$ of $\schc$ containing the zero object and all $\sigma$-semistable objects with phase $\phi$. One can also define for each interval $I \subset \br$ a subcategory $\schp(I)$ of $\schc$ as the extension-closed subcategory generated by all $\schp(\phi)$ with $\phi \in I$. According to \cite{Bri07}, the category $\schp(\phi)$ is abelian and $\schp(I)$ is quasi-abelian. Moreover, one can see $\schp((0, 1])= \scha$ by definition.

\subsection{The manifold of stability conditions}
A numerical stability condition is called \emph{locally finite} if there exists some $\epsilon>0$ such that for any $\phi\in\br$ the quasi-abelian subcategory $\schp(\phi-\epsilon,\phi+\epsilon)$ has finite length (for details, see \cite{Bri07}).

Suppose in addition that $\Num(\schc)$ has finite dimension, the set $\Stab(\schc)$ of all locally finite numerical stability conditions on $\schc$ admits a structure of complex manifold according to \cite{Bri07}. Henceforth, $\Num(\schc)$ is assumed to have finite dimension and the term \emph{stability condition} means a locally finite numerical stability condition.

\begin{proposition}[{{\cite[Section 9]{Bri08} and \cite[Proposition 2.3]{BM14}}}]\label{wall-and-chamber}
Given $v\in\Num(\schc)$, then there exists a locally finite set of walls (real codimension one submanifolds with boundary) in the complex manifold $\Stab(\schc)$, depending only on $v$, such that:
\begin{itemize}
	\item The sets of $\sigma$-semistable and $\sigma$-stable objects of class $v$ remain unchanged for any stability condition $\sigma$ in the same chamber;
	\item When $\sigma$ lies on a single wall in $\Stab(\schc)$, then there is a $\sigma$-semistable object that is unstable in one of the adjacent chambers, and semistable in the other adjacent chamber.
\end{itemize}
If $v$ is primitive in $\Num(\schc)$, then a stability condition $\sigma$ lies on a wall with respect to $v$ if and only if there exists a strictly $\sigma$-semistable object of class $v$. 
\end{proposition}

\begin{definition}
A stability condition $\sigma$ is called $v$\defi{-generic} for a given $v\in\Num(\schc)$ if it is not contained in any wall with respect to $v$.
\end{definition}

According to \cite[Lemma 8.2]{Bri07}, there are two mutually-commutative group actions on the complex manifold $\Stab(\schc)$ which are defined as follows.

The group $\Aut(\schc)$ acts on $\Stab(\schc)$ on the left, via 
$$\Phi. (\scha, Z):= (\Phi(\scha), Z \circ \Phi_*^{-1})  $$
where $\Phi_*$ denotes the induced automorphism of the group $\Num(\schc)$.

The universal cover 
$$
\grp  = \left \{ \tilde{g}=(M, f) ~  \middle\vert ~ 
\begin{array}{l} 
	M \in \GLp(2,\br), f \colon \br \rightarrow \br \text{ is an increasing function} \\ 
	\text{such that for all $\phi \in \br$ we have $f(\phi+1) = f(\phi) + 1$} \\
	\text{and $M \cdot e^{\pi \phi\sqrt{-1}} \in \br_{>0} \cdot e^{ \pi f(\phi)\sqrt{-1}}$} 
\end{array}  \right \}
$$
of the group $\GLp(2,\br)=\{M\in\GL(2,\br)\,|\,\det(M)>0\}$ acts on the right via 
$$
(\scha, Z) . (M, f):= (\schp((f(0), f(1)]), M^{-1} \circ Z). 
$$ 
where we use $\bc=\br\oplus\br\sqrt{-1}$ to validate the composite $M^{-1}\circ Z\colon \Num(\schc)\rightarrow\bc$.

\subsection{Stability conditions on K3 surfaces}
Let $S$ be a smooth projective K3 surface, then there exists a canonical polarized weight-two Hodge structure on the lattice $\tilde{H}(S,\bz)=H^0(S,\bz)\oplus H^2(S,\bz)\oplus H^4(S,\bz)$ whose polarization is given by the Mukai pairing $\langle(r_1,\alpha_1,s_1),(r_2,\alpha_2,s_2)\rangle_{\tilde{H}}=\alpha_1.\alpha_2-r_1.s_2-r_2.s_1$. The product $\alpha_1\cdot\alpha_2$ is the usual intersection product in $H^2(S,\bz)$. The Mukai vector 
$$v(E)=\ch(E).\sqrt{\td(S)}=(\rank(E),c_1(E),\rank(E)+c_1(E)^2/2-c_2(E))$$
of a complex $E\in\cate{D}^b(S)$ belongs to the algebraic part $\tilde{H}_{\textup{alg}}(S,\bz)$ of $\tilde{H}(S,\bz)$ and induces an isomorphism between $\tilde{H}_{\textup{alg}}(S,\bz)=H^0(S,\bz)\oplus \NS(S)\oplus H^4(S,\bz)$ and the numerical Grothendieck group $\Num(S)=\Num(\cate{D}^b(S))$.

An ample class $\omega\in \NS(S)_{\br}$ defines a notion of slope stability, called the $\mu_{\omega}$\emph{-stability}, by declaring the slope $\mu_{\omega}(\shf):=c_1(\shf)\cdot\omega/\rank(\shf)$ for a torsion-free sheaf $\shf$ on $S$. Truncating the HN filtrations at the number $\beta\cdot\omega$, a pair $(\omega,\beta)\in\NS(S)^2_{\br}$ with $\omega$ ample determines a torsion pair $(\scht_{\omega,\beta},\schf_{\omega,\beta})$ on $\cate{Coh}(S)$ such that
\begin{align*}
\scht_{\omega,\beta}&=\{\textup{torsion free part of }\shf\textup{ has }\mu_{\omega}\textup{-semistable HN factors of slope} >\beta\cdot\omega\}\\
\schf_{\omega,\beta}&=\{\shf \textup{ is torsion free and its }\mu_{\omega}\textup{-semistable HN factors have slope} \leq\beta\cdot\omega\}
\end{align*}
Here and henceforth, we adopt the convention that $\mu_{\omega}(\shf)=+\infty$ if $\shf$ is a torsion sheaf. This torsion pair gives a bounded $t$-structure on $\cate{D}^b(S)$ with heart 
$$\scha_{\omega,\beta}=\{E\in \cate{D}^b(S)\,|\, H^i(E)=0\textup{ for }i\notin\{-1,0\},H^{-1}(E)\in\schf_{\omega,\beta},H^0(E)\in\scht_{\omega,\beta}\}$$
according to \cite{HRS}. One recalls that an object $E$ in $\cate{D}^b(S)$ is called \emph{spherical} if it satisfies $\deriver\Hom(E,E)=\bc[0]\oplus\bc[-2]$. The group homomorphism
$$Z_{\omega,\beta}\colon \Num(S)\rightarrow\bc,\quad (r,\Delta,s)\mapsto \frac{1}{2}(2\beta\cdot\Delta-2s+r(\omega^2-\beta^2))+(\Delta-r\beta)\cdot\omega\sqrt{-1}$$
is a stability function on $\scha_{\omega,\beta}$ once $Z_{\omega,\beta}(\shf)\notin\br_{\leq0}$ for any spherical sheaf $\shf$ on $S$ due to \cite[Lemma 6.2]{Bri08}. Under this condition, the pair $\sigma_{\omega,\beta}=(\scha_{\omega,\beta},Z_{\omega,\beta})$ is a stability condition by \cite[Section 7]{Bri08}. Moreover, the skyscraper sheaf $\sho_p$ at any point $p\in S$ is $\sigma_{\omega,\beta}$-stable of phase $1$ for every $p\in S$.

A stability condition $\sigma$ on $S$ is called \emph{geometric} if the skyscraper sheaf $\sho_p$ is a $\sigma$-stable object of the same phase for every point $p\in S$. A geometric stability condition on $\cate{D}^b(S)$ has the form $\sigma_{\omega,\beta}.\tilde{g}$ for a unique determined stability condition $\sigma_{\omega,\beta}$ and a unique element $\tilde{g}\in\grp$ according to \cite[Proposition 10.3]{Bri08} and its argument. Also, there exists a connected component of $\Stab(\cate{D}^b(S))$ containing all geometric stability conditions which will be denoted by $\Stab^{\circ}(S)$.

\subsection{Autoequivalences on K3 surfaces and stability conditions}
According to \cite{Or97}, an autoequivalence of $\cate{D}^b(S)$ is isomorphic to a Fourier--Mukai transform.

\begin{definition}
Let $X$ and $Y$ be two smooth projective varieties and $K$ be an object of $\cate{D}^b(X\times Y)$, then the \defi{Fourier--Mukai transform with kernel $K$} is 
$$\Phi_K\colon \cate{D}^b(X)\rightarrow\cate{D}^b(Y),\quad E\mapsto \textsf{R}p(q^*E\otimes^{\textsf{L}} K)$$
where $q\colon X\times Y\rightarrow X$ and $p\colon X\times Y\rightarrow Y$ are the projections
\end{definition}

\begin{example}
Any spherical object $E$ in $\cate{D}^b(S)$ induces an autoequivalence $\funct{T}_{E}$ called the \emph{spherical twist} associated with $E$, whose Fourier--Mukai kernel is the mapping cone $\textsf{C}(p^*E^\vee\otimes^{\derivel} q^*E\rightarrow\sho_{\Delta})$ in $\cate{D}^b(S\times S)$. Here, $\sho_{\Delta}$ is the structure sheaf of the diagonal $\Delta\subset S\times S$ and $E^\vee$ is the derived dual of $E$. 
\end{example}

Similarly, one defines the \emph{cohomological Fourier--Mukai transform} $\Phi^{\coho}_K$ using the pushforwards and pullbacks in singular cohomology. It is proved in \cite{HMS09} that the assignment $\Phi_K\mapsto\Phi^{\coho}_K$ induces a surjection $\rho\colon\Aut(\cate{D}^b(S))\rightarrow \Aut^+(\tilde{H}(S,\bz))$ onto the group of isometries with positive orientation. The kernel $\Aut^{\circ}(\cate{D}^b(S))$ of the surjection $\rho$ is highly non-trivial due to the existence of spherical twists and it preserves $\Stab^{\circ}(S)$ when $S$ is generic (see \cite{BB17}). On the other hand, with the help of stability conditions, one can have a better understanding of $\Aut(\cate{D}^b(S))$.

\begin{proposition}[{{\cite[Proposition 1.4 and Corollary 1.5]{Hu16}}}]\label{K3-group-Aut(sigma)}
Consider a stability condition $\sigma=(\scha,Z)$ in $\Stab^{\circ}(S)$, then the following groups are finite
\begin{align*}
		\Aut(\cate{D}^b(S),\sigma)&:=\{\Phi\in\Aut(\cate{D}^b(S))\,|\,\Phi.\sigma=\sigma\}\\
		\Aut^+(\tilde{H}(S,\bz),Z)&:=\{\phi\in\Aut^+(\tilde{H}(S,\bz))\,|\,Z\circ\phi=Z\}
\end{align*}
and  $\rho$ induces an isomorphism between them.
\end{proposition}

Nevertheless, thanks to \cite{BB17}, a stronger statement and an inverse result have been recently obtained by Yu-Wei Fan and Kuan-Wen Lai.

\begin{proposition}[{{\cite[Theorem 4.13 and Lemma 5.2]{FL23}}}]\label{K3-group-realizing}
Consider a generic polarized K3 surface $S$ and a stability condition $\sigma\in\Stab(S)$, then $\Aut(\cate{D}^b(S),\sigma)$ is either trivial or isomorphic to $\bz/2\bz$. Conversely, any involution in $\Aut(\cate{D}^b(S))$ fixes a stability condition in $\Stab(S)$.
\end{proposition}

It is worth mentioning that $\Stab(S)=\Stab^{\circ}(S)$ for a generic polarized K3 surface $S$ according to \cite{BB17}.

\subsection{Moduli spaces of semistable objects}
An admissible subcategory $\schd$ inside $\cate{D}^b(X)$ of a smooth projective variety $X$ admits finite dimensional $\Num(\schc)$. Let $\sigma$ be a stability condition on $\schd$ and $v$ be a class in $\Num(\schd)$, one can define a moduli stack $\mathfrak{M}_{\sigma}(v)$ of $\sigma$-semistable objects with class $v$. A comprehensive reference is \cite{BLMNPS21}.

\begin{definition}
The moduli stack $\mathfrak{M}_{\sigma}(v)\colon(\cate{Sch}/\bc)^{\cate{op}}\rightarrow\cate{Gpds}$ is defined by
$$\mathfrak{M}_{\sigma}(v)(T)=\{E\in\cate{D}^b(X_T)\,|\,E|_{X_t}\in\schd\subset\cate{D}^b(X_t)\textup{ is }\sigma\textup{-semistable for any }t\in T\}$$
for any locally finitely generated complex scheme $T$. Similarly, one defines $\mathfrak{M}^{\stable}_{\sigma}(v)$.
\end{definition}

Given a class $v\in\Num(S)$ and a $v$-generic stability condition $\sigma$ in $\Stab^{\circ}(S)$ for a K3 surface $S$, then \cite{BM14} asserts that the moduli stack $\mathfrak{M}_{\sigma}(v)$ admits a projective coarse moduli space $M_{\sigma}(v)$ whose points are $S$-equivalent classes of $\sigma$-semistable objects in $\schp(-1,1]$ with Mukai vector $v$. The moduli space $M_{\sigma}(v)$ does not change if we move $\sigma$ within a chamber with respect to $v$  (see Proposition \ref{wall-and-chamber}) and behaviors of $M_{\sigma}(v)$ with $\sigma$ moving across a wall are studied intensively in \cite{BM14,BM14-MMP,MZ16}. Here we only recall some necessary definitions for this article.

Suppose that $W\subset\Stab^{\circ}(S)$ is a wall with respect to a given $v\in \Num(\schd)$ with $M_{\sigma_0}(v)\neq0$ and $\sigma_0\in W$ is not in any other wall, then $W$ is called a \emph{totally semistable wall} if the stable locus $M^{\stable}_{\sigma_0}(v)$ is empty.

Suppose in addition that $\sigma_{\pm}$ are $v$-generic stability conditions in two different adjoining chambers of the wall $W$. Then \cite[Proposition 5.2]{MZ16} asserts a birational morphism $\pi_{\pm}\colon M_{\sigma_{\pm}}(v)\rightarrow\overline{M}_{\pm}$ which contracts curves in $M_{\sigma_{\pm}}(v)$ parameterizing $S$-equivalent objects under $\sigma_0$, where $M_{\pm}$ is the image of $\pi_{\pm}$ in $M_{\sigma_0}(v)$. Then the wall $W$ is called a \emph{fake wall} if $\pi_{\pm}$ are isomorphisms, a \emph{flopping wall} if $\pi_{\pm}$ are small contraction, and a \emph{divisorial wall} if $\pi_{\pm}$ are divisorial contractions.

\subsection{Mukai isomorphism} In this subsection, we will restrict ourselves on the moduli space $M_{\sigma}(v_n)$ for the Mukai vector $v_n=(1,0,1-n)$. It is known that for any ample divisor $H$ on a K3 surface, any Gieseker $H$-semistable sheaf $\shg$ with $v(\shg)=v_n$ is isomorphic to the ideal sheaf $\shi_Z$ of a zero-dimensional subscheme of $S$ with length $n$ and the moduli space $M_H(v_n)$ of Gieseker $H$-semistable sheaves is isomorphic to $S^{[n]}$. By \cite[Proposition 14.2]{Bri08} and \cite[Section 6]{Yo01}, one also has

\begin{proposition}\label{K3-large-volume-limits}
Consider an ample class $H$ on $S$ and real classes $\omega,\beta\in\NS(S)_{\br}$ such that $\omega= rH$ for some $r>0$ and $\sigma_{\omega,\beta}$ is a stability condition, then there exists a positive number $t_0$ such that, for any $t\geq t_0$, an object $E$ of $\cate{D}^b(S)$ with $v(E)=v_n$ is $\sigma_{t\omega,\beta}$-semistable if and only if $E$ is isomorphic to a shift of a Gieseker $H$-semistable sheaf. Moreover, one has $M_{\sigma_{t\omega,\beta}}(v_n)=M_H(v_n)\cong S^{[n]}$ for $t\gg0$.
\end{proposition}

In general, one has the following application of \cite[Corollary 6.9]{BM14}.

\begin{proposition}\label{K3-moduli-stable-hyperkaehler}
Consider a $v_n$-generic stability condition $\sigma$ in $\Stab^{\circ}(S)$, then $M_{\sigma}(v_n)=M^{\stable}_{\sigma}(v_n)$ is a hyperkähler manifold deformation equivalent to $S^{[n]}$.
\end{proposition}

In this case, the space $H^2(M_{\sigma}(v_n),\bz)$ is torsion-free and admits a polarized weight-two Hodge structure according to \cite{Bea83}. The orthogonal space $v_n^{\perp}:=\{w\in \tilde{H}(S,\bz)\,|\,\langle v_n,w\rangle_{\tilde{H}}=0\}$ also inherits a Hodge structure from $\tilde{H}(S,\bz)$ polarized by the Mukai pairing. It turns out that these two Hodge structures are isometric.

\begin{proposition}[{{\cite[Theorem 6.10]{BM14}}}]
Consider a $v_n$-generic stability condition $\sigma\in\Stab^{\circ}(S)$ such that $M_{\sigma}(v)$, then there exists a canonical Hodge isometry
$$\theta_{\sigma,n}\colon v_n^{\perp}\rightarrow H^2(M_{\sigma}(v_n),\bz)$$
which only depends on the chamber of $\sigma$ with respect to $v_n$.
\end{proposition}

This isometry is called the \emph{Mukai isomorphism} and is firstly introduced by Mukai in \cite{Mu87} for moduli space of semistable sheaves.

\begin{example}\label{K3-hilbert-realization}
Under the Mukai isomorphism $\theta_n\colon v_n^{\perp}\rightarrow H^2(S^{[n]},\bz)$, the algebraic part of $H^2(S^{[n]},\bz)$ always contains $\theta_n(1,0,n-1)$ and $\theta_n(0,H,0)$.
\end{example}

Moreover, the Mukai isomorphism is compatible with cohomology and we provide here a special case of \cite[Proposition 3.5]{Ou18}. The original statement in \cite{Ou18} is only for isomorphisms, but the argument also works for birational maps.

\begin{proposition}\label{K3-invariant-lattices-compatible}
Consider an autoequivalence $\Phi_K$ of $\cate{D}^b(S)$ with $\Phi_K^{\coho}(v_n)=v_n$ and two $v_n$-generic stability condition $\sigma,\tau\in\Stab^{\circ}(S)$. Suppose that $\Phi_K$ induces a biregular isomorphism $\phi\colon M_{\sigma}(v_n)\rightarrow M_{\tau}(v_n)$ (resp.\ a birational map $\phi\colon M_{\sigma}(v_n)\dashrightarrow M_{\tau}(v_n)$) by the mapping $[E]\mapsto[\Phi(E)]$, then
\begin{displaymath}
	\xymatrix{
		v_n^{\perp}\ar[d]_{\theta_{\sigma,n}}\ar[rr]^{\Phi_K^{\coho}}&&v_n^{\perp}\ar[d]_{\theta_{\tau,n}}\\
		H^2(M_{\sigma}(v_n),\bz)\ar[rr]^{\phi_*}&&H^2(M_{\tau}(v_n),\bz)
	}
\end{displaymath}
is a commutative diagram, where $\phi_*$ is the Hodge isometry induced by $\phi$.
\end{proposition}

\section{Proof of Theorem \ref{1-main-result} and Proposition \ref{1-non-main-result}}
\subsection{The necessary part and valid automorphisms}
At first, we introduce the definition and show the necessary part of Theorem \ref{1-main-result} in a general context. 

\begin{definition}\label{main-definition}
A biregular (resp.\ birational) automorphism $f$ on $S^{[n]}$ for some $n\geq 2$ is called \defi{derived-natural}, if there exists $\Phi\in\cate{Aut}(\cate{D}^b(S))$ such that $\Phi(\shi_Z)\cong\shi_{f(Z)}$ for any $Z\in S^{[n]}$ (resp.\ for any $Z$ in an open regular locus of $f$).
\end{definition}

One has an injection $H^2(S^{[n]},\bz)\hookrightarrow H^2(S^{[n]},\bz)^\vee$ by sending $\alpha$ to $\langle\alpha,-\rangle_{\tilde{H}}$ and its cokernel is the discriminant group $A_{H^2(S^{[n]},\bz)}$ of $H^2(S^{[n]},\bz)$. A Hodge isometry of $H^2(S^{[2]},\bz)$ induces an automorphism on $A_{H^2(S^{[n]},\bz)}$ (see e.g.\ \cite[14.2.2]{HuK3}).

\begin{proposition}\label{necessary condition}
Consider a K3 surface $S$ and a derived-natural automorphism $f$ of $S^{[n]}$, then the Hodge isometry $f_*$ induces the identity on $A_{H^2(S^{[n]},\bz)}$.
	
\begin{proof}
Let $\Phi_K$ be an autoequivalence that induces $f$, then by virtue of Proposition \ref{K3-invariant-lattices-compatible} the Hodge isometry $f_*$ on $v_n^{\perp}\subset \tilde{H}(S,\bz)$ extends to the Hodge isometry $\tilde{H}(S,\bz)$ by claiming $\Phi_K^{\coho}(v_n)=v_n$. So the induced automorphism on $A_{H^2(S^{[n]},\bz)}$ by $f_*$ is the identity by standard facts in lattice theory (see e.g.\ \cite[Proposition 14.2.6]{HuK3}).
\end{proof}
\end{proposition}

Now we come to generic polarized K3 surfaces. The biregular automorphisms on Hilbert schemes of $n$ points on a generic K3 surface are worked out in \cite{Ca19} and the birational ones are studied in \cite{BC22}. Here we only cite the relevant parts.

\begin{proposition}\label{involution-non-natural-derived-natural}
Consider a generic polarized K3 surface $(S,H)$ of degree $2t$ with $t\geq2$ and an integer $n\geq 2$, then $S^{[n]}$ admits a non-trivial birational involution $\phi$ such that $\phi_*$ induces the identity on $A_{H^2(S^{[n]},\bz)}$ if and only if the integer $t(n-1)$ is not a square and the negative Pell's equation
\begin{equation}\label{Pell-derived-natural}
	X^2-t(n-1)Y^2=-1
\end{equation}
admits solutions. In this case, the isometry $\phi_*$ of the sublattice $\NS(S^{[n]})$ is represented by the matrix
$$\begin{pmatrix}
	2a^2+1&-2(n-1)ab\\2tab&-2a^2-1
\end{pmatrix}$$
under the basis $\{\theta_n(0,-H,0),\theta_n(1,0,n-1)\}$, where $(a,b)$ is the positive solution of Equation (\ref{Pell-derived-natural}) with minimal $a>0$. 
\end{proposition}

\begin{remark}\label{remarkn}
Under the same conditions, the equation $(n-1)X^2-tY^2=1$ does not have a solution for $n\neq2$ by \cite[Lemma A.1]{BC22}. It means that the condition \cite[Theorem 1.3(ii)]{Ca19} is always true in this case. Moreover, $(2a^2+1,2ab)$ is the minimal positive solution of the equation $X^2-t(n-1)Y^2=1$.
\end{remark}

\begin{remark}\label{involution-regular}
By the proof of \cite[Theorem 6.4]{Ca19}, the birational involution $\phi$ is biregular if there are no flopping walls inside the movable cone of $S^{[n]}$. There is a numerical description of this condition which is stated in \cite[Theorem 1.3(iii)]{Ca19}.
\end{remark}

The degree $2$ case is different, because such a K3 surface $S$ is a branched double cover over $\bp^2$ and has a covering involution $\tau$. The group $\Aut(S^{[n]})$ is described as following by the arguments for \cite[Theorem 1.1 and Proposition 2.2]{BC22}.

\begin{proposition}\label{involution-degree2}
Consider a generic polarized K3 surface $(S,H)$ of degree $2$, then
$$\Aut(S^{[n]})=\{\phi^{[n]},\id\}$$
where $\phi^{[n]}$ is the involution induced by $\tau$. Moreover, any other possible non-natural birational automorphism of $S^{[n]}$ acts on $A_{H^2(S^{[n]},\bz)}$ by $-1$.
\end{proposition}

\subsection{The existence of stability conditions and involutions}
It remains to consider the automorphism $\phi$ stated in Proposition \ref{involution-non-natural-derived-natural}. In this case, the isometry $\phi_*$ can be extended to a Hodge isometry $\tau$ on $\tilde{H}(S,\bz)$ via $\theta_n$ by claiming that $\tau(v_n)=(v_n)$ (see e.g.\ \cite[Proposition 14.2.6]{Hu}). Using the matrix in Proposition \ref{involution-non-natural-derived-natural}, one computes directly that the restriction of $\tau$ on $\Num(S)$ is given by the matrix
$$\begin{pmatrix}
	-a^2&-2tab&-tb^2\\
	(n-1)ab&2a^2+1&ab\\
	-(n-1)^2tb^2&-2t(n-1)ab&-a^2
\end{pmatrix}$$
with respect to the basis $(1,0,0),(0,H,0),(0,0,1)$.

\begin{proposition}\label{stabilityconditionpreexist}
Under the setups in Proposition \ref{involution-non-natural-derived-natural} and above, the equation
$$Z_{\omega,\beta}(\tau(r,mH,s))=Z_{\omega,\beta}(r,mH,s)$$ 
holds for any $(r,mH,s)\in\Num(S)$ if and only if $\omega=\frac{1}{tb}H$ and $\beta=-\frac{a}{tb}H$.
	
\begin{proof}
Comparing the imaginary parts, one obtains an equation
$$((n-1)br+2am+bs)(aH\cdot\omega+tb\omega\cdot\beta)=0$$
which directly implies $aH\cdot\omega+tb\omega\cdot\beta=0$ or equivalently $aH+tb\beta=0$.
		
Comparing the real parts, one obtains an equation
$$
((n-1)br+2am+bs)\left(\frac{tb}{2}\omega^2+\frac{a^2}{b}-t(n-1)b\right)=0
$$
so that $tb^2\omega^2=2t(n-1)b^2-2a^2=2$ and hence $tb\omega=H$.
\end{proof}
\end{proposition}

\begin{proposition}\label{stabilityconditionexist}
Under the above setups, the pair $(\omega_0,\beta_0)=(\frac{1}{tb} H,-\frac{a}{tb} H)$ gives a stability condition $\sigma_{\omega_0,\beta_0}$.

\begin{proof}
According to Section 2.3, it suffices to check $\Re(Z_{\omega_0,\beta_0}(r,mH,s))>0$ for any $v=(r,mH,s)$ satisfying $r>0,\langle v,v\rangle_{\tilde{H}}=-2,\Im(Z_{\omega_0,\beta_0}(v))=0$. One has
$$Z_{\omega_0,\beta_0}(r,mH,s)=\left(-\frac{2a}{b}m-s+(\frac{1}{tb^2}-\frac{a^2}{tb^2})r\right)+\left(\frac{2a}{tb^2}r+ \frac{2}{b}m\right)\sqrt{-1}$$
from which one reduces to show $ar+tmb=0,tm^2=rs-1,r>0\Rightarrow (n-1)r>s$.
	
At first, one notices that
$$tm=-\frac{a}{b}r<0\Rightarrow rs-1=tm^2\leq \frac{a^2}{tb^2}r^2\Leftrightarrow \frac{a^2}{a^2+1}(n-1)r+\frac{1}{r}\geq s$$
from which one obtains $(n-1)r\geq s$. In fact, one has
$$1>\frac{1}{r}-\frac{1}{a^2+1}(n-1)r\geq s-(n-1)r$$

It remains to show that $(n-1)r=s$ is not possible. Otherwise, one immediately has $(n-1)r^2-tm^2=1$. It contradicts Remark \ref{remarkn} when $n\neq2$. Suppose that $n=2$, one has $tm^2(a^2+1)=a^2r^2$ so that $a^2r^2=ta^2m^2+a^2=ta^2m^2+tm^2$ and then $r^2-1=tm^2=a^2$. It can never happen as $r,a$ are positive integers.
\end{proof}
\end{proposition}

\begin{remark}
Similarly, one can show that $\sigma_{\lambda\omega_0,\beta_0}$ is a stability condition for any
real number $\lambda\geq1$. Moreover, the pair $(\omega_0,\beta_0)$ might still give a stability condition for non-generic polarized K3 surfaces (see e.g.\ \cite{Liu24}).
\end{remark}

By virtue of Proposition \ref{K3-group-Aut(sigma)} and \ref{stabilityconditionpreexist}, there exists an involution $\Phi$ on $\cate{D}^b(S)$ realizing the Hodge isometry $\tau$ and fixing $\sigma_{\omega_0,\beta_0}$. We will elaborate on $\Phi$ later.

\subsection{Identification of the moduli spaces via wall-crossing}
We will next show that, for any given $\lambda\geq 1$, the stability condition $\sigma_{\lambda\omega_0,\beta_0}$ is either $v_n$-generic or lies in a flopping wall with respect to $v_n$. To achieve this, we need to exclude the existence of all the three types divisorial walls: the Brill--Noether walls, the Hilbert--Chow walls, and the Li--Gieseker--Uhlenbeck walls, as well as existence of the totally semistable walls. One proceeds by applying \cite[Theorem 5.7]{BM14-MMP}.

A potential wall for $\sigma_{\lambda\omega_0,\beta_0}$ with respect to $v_n$ is associated with a lattice $\mathbb{H}$ such that any $v=(r,mH,s)\in\mathbb{H}$ satisfies $\Im \left(Z_{\sigma_{\lambda\omega_0,\beta_0}}(v_n)/Z_{\sigma_{\lambda\omega_0,\beta_0}}(v)\right)=0$. It is also
$$\left(-\frac{2a}{b}m-s+(\frac{\lambda^2}{tb^2}-\frac{a^2}{tb^2})r\right)\frac{2a}{tb^2}=\left(\frac{2a}{tb^2}r+ \frac{2}{b}m\right)\left(n-1+\frac{\lambda^2}{tb^2}-\frac{a^2}{tb^2}\right)$$
or equivalently $(2(n-1)b^2t+\lambda^2-1)m+ab((n-1)r+s)=0$. To save notation, we define $\lambda_0\geq0$ by $2(n-1)b\lambda_0=\lambda^2-1$. Then 
\begin{equation}\label{vectors-in-lattice}
(n-1)r+s=-(bt+\lambda_0)\frac{2(n-1)m}{a}
\end{equation}

\begin{proposition}\label{WC-BN-wall}
There exists no vector $w\in\mathbb{H}$ with $\langle w,w\rangle_{\tilde{H}}=-2$ and $\langle w,v_n\rangle_{\tilde{H}}=0$, so that the stability condition $\sigma_{\lambda\omega_0,\beta_0}$ does not lie on a Brill--Noether wall.

\begin{proof}
Suppose otherwise that there exists such a vector $w=(r,mH,s)$. Then one has $(n-1)r=s$ and $rs-1=tm^2$. It follows that $(n-1)r^2-tm^2=1$, which contradicts to Remark \ref{remarkn} when $n\neq 2$. Suppose in addition that $n=2$, then by Equation \ref{vectors-in-lattice} one has
$$r=s=-(bt+\lambda_0)\frac{m}{a}$$
which, combing with the equalities $r^2-tm^2=1$ and $a^2-tb^2=-1$, implies
$$tm^2+2bt\lambda_0m^2+\lambda_0^2m^2=(tb^2-a^2)tm^2+2bt\lambda_0m^2+\lambda_0^2m^2=a^2$$
On the other hand, $r^2-tm^2=1$ and Remark \ref{remarkn} ensure that $|r|\geq 2a^2+1,|m|\geq ab$. So one can immediately see a contradiction from $tm^2\geq ta^2b^2>a^2$.
\end{proof}
\end{proposition}

\begin{proposition}\label{WC-HC-wall}
There exists no vector $w\in\mathbb{H}$ with $\langle w,w\rangle_{\tilde{H}}=0$ and $\langle w,v_n\rangle_{\tilde{H}}=1$, so that $\sigma_{\lambda\omega_0,\beta_0}$ does not lie on a Hilbert--Chow wall. 

\begin{proof}
Suppose otherwise that there exists such a vector $w=(r,mH,s)$. Then 
$$(2(n-1)r-1)^2-t(n-1)(2m)^2=1$$ 
from which it follows by Remark \ref{remarkn} that
$$|(n-1)r+s|=|2(n-1)r-1|\geq2a^2+1 \quad\textup{and}\quad |2m|\geq 2ab$$
On the other hand, one has
$$r=-\frac{(bt+\lambda_0)m}{a}+\frac{1}{2(n-1)}\quad\textup{and}\quad s=-\frac{(bt+\lambda_0)m(n-1)}{a}-\frac{1}{2}$$
so that $4(t^2b^2+2\lambda_0tb+\lambda_0^2)m^2(n-1)^2-4t(n-1)m^2a^2-a^2=0$ and
$$\left|(bt+\lambda_0)\frac{2(n-1)m}{a}\right|\geq 2a^2+1\quad\textup{or equivalently}\quad \lambda_0\geq\frac{(2a^2+1)a}{2(n-1)|m|}-tb$$

However, the function $f(x)=4(t^2b^2+2tbx+x^2)m^2(n-1)^2-4t(n-1)m^2a^2-a^2$ has a zero in that range if and only if $|2(n-1)r-1|=2a^2+1$ and $|m|=ab$. Therefore, one has
$$\lambda_0=\frac{2a^2+1}{2(n-1)b}-tb=\frac{2a^2-2(n-1)tb^2+1}{2(n-1)b}=-\frac{1}{2(n-1)b}<0$$
causing a contradiction to $\lambda_0\geq0$
\end{proof}
\end{proposition}

\begin{proposition}\label{WC-LHU-wall}
There exists no vector $w\in\mathbb{H}$ with $\langle w,w\rangle_{\tilde{H}}=-2$ and $\langle w,v_n\rangle_{\tilde{H}}=2$, so that $\sigma_{\lambda\omega_0,\beta_0}$ does not lie on a Li--Gieseker--Uhlenbeck wall. 

\begin{proof}
Suppose otherwise that there exists such a vector $w=(r,mH,s)$. Then 
$$((n-1)r-1)^2-t(n-1)m^2=1$$ 
from which it follows by Remark \ref{remarkn} that
$$|(n-1)r-1|\geq2a^2+1 \quad\textup{and}\quad |m|\geq 2ab$$
On the other hand, one has
$$r=-\frac{(bt+\lambda_0)m(n-1)}{a}+\frac{1}{(n-1)}\quad\textup{and}\quad s=-\frac{(bt+\lambda_0)m(n-1)}{a}-1$$
so that $(t^2b^2+2\lambda_0tb+\lambda_0^2)m^2(n-1)^2-t(n-1)m^2a^2-a^2=0$ or equivalently
$$(n-1)tm^2+(2tb+\lambda_0)\lambda_0m^2(n-1)^2=a^2$$
which is impossible as $|m|\geq 2ab,2tb>0,\lambda_0\geq0$.
\end{proof}
\end{proposition}

To tackle the last possible non-flopping wall, we need to introduce the \emph{effective vectors}. A vector $w$ in $\Num(S)$ associated with $\sigma_{\lambda\omega_0,\beta_0}$ is called \emph{effective} if it satisfies $\langle w,w\rangle_{\tilde{H}}\geq-2$ and $\Re(Z_{\sigma_{\lambda\omega_0,\beta_0}}(w)/Z_{\sigma_{\lambda\omega_0,\beta_0}}(v_n))>0$.

\begin{proposition}\label{WC-ss}
There exists no effective vector $w\in\mathbb{H}$ with $\langle w,w\rangle_{\tilde{H}}=-2$ and $\langle w,v_n\rangle_{\tilde{H}}<0$, so that $\sigma_{\lambda\omega_0,\beta_0}$ does not lie on a totally semistable wall.
	
\begin{proof}
Suppose otherwise that there exists such a vector $w=(r,mH,s)$. Then one has $(n-1)rs=(n-1)(tm^2+1)$ and $s>(n-1)r$. Therefore, one obtains
$$(n-1)r=-(bt+\lambda_0)\frac{(n-1)m}{a}-\sqrt{\frac{(tb+\lambda_0)^2(n-1)^2m^2}{a^2}-(n-1)(tm^2+1)}$$ $$s=-(bt+\lambda_0)\frac{(n-1)m}{a}+\sqrt{\frac{(tb+\lambda_0)^2(n-1)^2m^2}{a^2}-(n-1)(tm^2+1)}$$
by also considering the Equation (\ref{vectors-in-lattice}).

The spherical vector $w$ is also effective, so one has
\begin{equation*}
\begin{split}
\left(-\frac{2a}{b}m-s+\frac{2(n-1)b\lambda_0+2-tb^2}{tb^2}r\right)&\frac{2(n-1)b\lambda_0+2}{tb^2}+\\
\frac{2a}{tb^2}&\left(\frac{2a}{tb^2}r+ \frac{2}{b}m\right)(2(n-1)b\lambda_0+1)>0\quad
\end{split}
\end{equation*}
which is also
\begin{equation*}
	\begin{split}
(4n-4)a^2b\lambda_0r+2a^2r+(2n-2)&ab^2mt\lambda_0-((n-1)b\lambda_0+1)b^2t(r+s)\\
&+2(n-1)^2b^2\lambda_0^2r+(4n-4)b\lambda_0r+2r>0
	\end{split}
\end{equation*}

One can then obtain, after plugging into $r,s$ and $a^2=(n-1)tb^2-1$, that $f(\lambda_0)((n-1)m\lambda_0+\sqrt{(tb+\lambda_0)^2(n-1)^2m^2-(n-1)a^2(tm^2+1)}<0$ where $f(\lambda_0)=2(n-1)^2b^2\lambda_0^2+(4(n-1)^2b^3t+(n-2)(bn-b+1))\lambda_0+2(n-1)b^2t>0$. Therefore
$$\sqrt{(tb+\lambda_0)^2(n-1)^2m^2-(n-1)a^2(tm^2+1)}<-m(n-1)\lambda_0$$ 
It means that $m<0$ and $(tb+\lambda_0)^2(n-1)m^2-a^2(tm^2+1)<(n-1)m^2\lambda_0^2$, which is also $(1+2(n-1)b\lambda_0)tm^2<a^2$. So one obtains a bound
$$m^2<\frac{a^2}{t+2(n-1)b\lambda_0t}.$$
On the other hand, to have valid integers $r$ and $s$, one needs 
$$s-(n-1)r=2\sqrt{\frac{(tb+\lambda_0)^2(n-1)^2m^2}{a^2}-(n-1)(tm^2+1)}\geq 1$$
which follows another bound
$$m^2\geq\frac{4n-3}{4(n-1)}\cdot\frac{a^2}{t+2(n-1)tb\lambda_0+(n-1)\lambda^2_0}.$$

Combining the two bounds for $m^2$, one obtains $\lambda_0>(bt+\sqrt{b^2t^2+4t})/4(n-1)$. Taking it back to the first bound for $m^2$, one has
$$0<m^2<\frac{4tb^2-4}{4t+2b^2t^2+2bt\sqrt{b^2t^2+4t}}<1$$ 
which contradicts the fact that $m$ is an integer.		
\end{proof}
\end{proposition}

The previous four propositions in this subsection ensures that $\sigma_{\lambda\omega_0,\beta_0}$ is either generic or lies on a non-totally-semistable flopping wall with respect to $v_n$. Such a flopping wall induces a birational map $M_{\sigma_-}(v_n)\dashrightarrow M_{\sigma_+}(v_n)$ for $v_n$-generic stability condition $\sigma_{\pm}$ in adjoining chambers by contracting some objects according to \cite[Theorem 1.4 (b)]{BM14}, so one obtains a birational map $M_{\sigma_{\omega_0,\beta_0}}(v_n)\dashrightarrow M_H(v_n)$ by Proposition \ref{K3-large-volume-limits}. Since $\Phi$ induces an involution on $M_{\sigma_{\omega_0,\beta_0}}(v_n)$, it also induces a birational involution on $M_H(v_n)$. This involution is precisely the $f$ we start with by virtue of Proposition \ref{K3-invariant-lattices-compatible}. So we have finished the birational part of Theorem \ref{1-main-result}. The biregular part needs the condition \cite[Theorem 1.3(iii)]{Ca19}.

\begin{proposition}
Suppose that $f$ is biregular, then $M_{\sigma_{\omega_0,\beta_0}}(v_n)=M_H(v_n)$.

\begin{proof}
Otherwise, there exists a vector $w=(r,mH,s)$ in $\mathbb{H}$ with $\langle w,v_n\rangle_{\tilde{H}}=k$ and $\langle w,w\rangle_{\tilde{H}}=2p$ such that the integers $p,k$ fall into one of the cases in \cite[Theorem 1.3(iii)]{Ca19} (which uses the notation $\rho,\alpha$ for integers). Since $(2(n-1)r-k)^2-4t(n-1)m^2=k^2-4p(n-1)$, one has
$$\frac{|m|}{|2(n-1)r-k|}\geq\frac{ab}{2a^2+1}$$
which follows directly that $\lambda_0<-1/2(n-1)b<0$ by applying Equation (\ref{vectors-in-lattice}).
\end{proof}
\end{proposition}

\subsection{More on inducing involutions} Here we prove Proposition \ref{1-non-main-result} and present some examples. An important ingredient is the following result.

\begin{theorem}\cite{BB17}
Consider a generic K3 surface $S$, then the kernel of the map
$$\rho\colon \Aut(\cate{D}^b(S))\rightarrow\Aut^+(\tilde{H}(S,\bz))$$
preserves the connected component $\Stab^{\circ}(S)$.
\end{theorem} 

This is a special case of the famous conjecture \cite[Conjecture 1.2]{Bri08}. The argument we will see subsequently in fact works whenever this conjecture is true.

\begin{proof}[The Proof of Proposition \ref{1-non-main-result}]
Choose two autoequivalences $\Phi$ and $\Psi$ inducing the same $f$, then the composite $\Phi\circ\Psi^{-1}$ should induce the identity morphism on $M_H(v_n)$. In particular, the Hodge isometry $(\Phi\circ\Psi^{-1})_*$ on $\tilde{H}(S,\bz)$ is the identity according to Proposition \ref{K3-invariant-lattices-compatible}. So $\Phi\circ\Psi^{-1}$ is in $\ker(\rho)$ which preserves $\Stab^{\circ}(S)$ according to \cite{BB17}. Hence, by \cite[Theorem 13.3]{Bri08}, it takes $\sigma$ to $\sigma.\tilde{g}$ for any $\sigma\in\Stab^{\circ}(S)$ and some $\tilde{g}\in\grp$ depending on $\sigma$. But the fact that $\Phi\circ\Psi^{-1}$ fixes the moduli space $M_H(v_n)$ forces $\tilde{g}$ to be trivial. Thus $\Phi\circ\Psi^{-1}$ is isomorphic to  $\funct{Id}$.
\end{proof}

Now we give a description of the unique autoequivalence by appealing to the construction used in Proposition \ref{K3-group-Aut(sigma)}.

\begin{proposition}
Consider a derived-natural involution in Proposition \ref{involution-non-natural-derived-natural}, then the autoequivalence $\Phi$ is isomorphic to the Fourier--Mukai transform $\Phi_{\shu}[1]$, where $\shu$ on $S\times M$ is a universal family for the fine moduli space $M= M_H(tb^2,-abH,a^2)$.
	
\begin{proof}
Using the argument of \cite[Corollary 10.12]{Hu}, there exists a universal family $\shu$ for the fine moduli space $M_H(tb^2,-abH,a^2)$ such that the Fourier--Mukai transform $\Psi:=\Phi_{\shu}[1]$ induces the Hodge isometry $\tau$ whose restriction on $H^2(S^{[n]},\bz)$ is $\phi_*$. Here we need to recall the fact that $\Aut(S)$ is trivial and there are no rational smooth curves on $S$. The autoequivalence $\Psi$ preserves $\Stab^{\circ}(S)$ due to \cite[Proposition 4.2]{Hu08}. Choose the stability condition $\sigma_{\omega_0,\beta_0}$ from Proposition \ref{stabilityconditionexist}, then one has $\Psi.\sigma_{\omega_0,\beta_0}.=\sigma_{\omega,\beta}.\tilde{g}$ for some stability condition $\sigma_{\omega,\beta}$ and $\tilde{g}\in\grp$. However, by inspecting the central charge, one has $\omega=\omega_0,\beta=\beta_0$ and $\tilde{g}=(I_2,f)$. It suffices to check that it preserves the heart of the stability condition. In other words, the shift $\she[1]$ of any stable vector bundle $\she$ in $M_H(tb^2,-abH,a^2)$ is in the heart $\scha_{\omega_0,\beta_0}$. It is true because $\mu_{\omega_0}(\she)=\omega_0.\beta_0$ by assumption and then $\she\in\schf_{\omega_0,\beta_0}$. So one has proved that $\Psi.\sigma_{\omega_0,\beta_0}=\sigma_{\omega_0,\beta_0}$ and can conclude by Proposition \ref{K3-group-Aut(sigma)}.
\end{proof}
\end{proposition}

Unfortunately, this description is not helpful to describe the derived-natural automorphisms point-wisely. A better description is only known in limited cases and the following examples are computed in \cite{Liu24} and are derived from \cite{KP17}.

\begin{example}\label{Example1}
Let $(S,H)$ be a generic K3 surface of degree $4$, then the unique biregular involution on $S^{[2]}$, also known as the Beauville involution \cite{Bea83b}, is derived-natural and is induced by $\funct{T}_{\sho_S}\circ(-\otimes\sho_S(H))\circ\funct{T}_{\sho_S}\circ(-\otimes\sho_S(H))[-1]$.
\end{example}

\begin{example}\label{Example2}
Let $(S,H)$ be a generic K3 surface of degree $10$, then the birational involution on $S^{[2]}$ in \cite[Section 4.3]{O'Grady05} is derived-natural. The corresponding autoequivalence is $\funct{T}_{\sho_S}\circ\funct{T}_{\shu^\vee_S}\circ(-\otimes\sho_S(H))[-1]$, where $\shu_S$ is the unique element in $M_H(2,-H,3)$ (the uniqueness can be seen from \cite[Lemma 7.1]{BM14}). Also, the birational involution on $S^{[3]}$ described in \cite[Example 4.12]{De22} is derived-natural and the corresponding autoequivalence is $\funct{T}_{\shu_S}\circ\funct{T}_{\sho_S}\circ(-\otimes\sho_S(H))[-1]$.
\end{example}

The autoequivalences in Example \ref{Example2} are conjugate. Such a phenomenon is studied thoroughly in \cite{FL23} and we cite the following special case. 

\begin{theorem}[{{\cite[Theorem 5.6]{FL23}}}]
Consider a generic polarized K3 surface $S$ of degree $4$ or $10$, then the non-trivial maximal finite subgroups of $\Aut(\cate{D}^b(S))$ are isomorphic to $\bz/2\bz$ and conjugate to each other.
\footnote{Special thanks to Kuan-Wen Lai for explaining this result to me.}
\end{theorem}

It implies that the inducing autoequivalence for a derived-natural automorphism on $S^{[n]}$ for a generic K3 surface $S$ of degree $4$ or $10$ can be described using the constructions in \cite{KP17,PPZ23} and computations in \cite{Liu24}. Here are two first examples.

\begin{example}\label{Example-degree4}
Let $(S,H)$ be a generic K3 surface of degree $4$, then the birational involution on $S^{[6]}$ described in \cite[Example 6.1]{BC22} point-wisely is derived-natural and induced by $\funct{T}_{\sho_S(-H)}\circ(-\otimes\sho_S(H))\circ\funct{T}_{\sho_S(-H)}\circ(-\otimes\sho_S(H))[-1]$.
\end{example}

\begin{example}\label{Example-degree10}
Let $(S,H)$ be a generic K3 surface of degree $10$, then there exists a unique derived-natural birational involution on $S^{[11]}$ by Proposition \ref{involution-non-natural-derived-natural}. One can obtain that the involution is induced by $\funct{T}_{\sho_S(-H)}\circ\funct{T}_{\shu^\vee_S(-H)}\circ(-\otimes\sho_S(H))[-1]$. Similarly, the unique birational involution on $S^{[14]}$ is derived-natural and is induced by the autoequivalence $\funct{T}_{\shu_S(-H)}\circ\funct{T}_{\sho_S(-H)}\circ(-\otimes\sho_S(H))[-1]$. These two involutions have no geometric descriptions in the literature.
\end{example}

\begin{remark}
Unfortunately, a similar construction does not work for generic K3 surfaces of higher degree according to \cite[Appendix A]{Liu24}. 
\end{remark}

To conclude, we comment that the idea of derived-natural automorphisms can be adjusted to other types of hyperkähler manifolds. In other words, some automorphisms of hyperkähler manifolds can be realized by autoequivalences on certain categories similar to $\cate{D}^b(S)$. It is a generalization to the notion of induced automorphisms studied in \cite{CKKM19,MW15} via lattices. Several examples of such automorphisms on $K3^{[n]}$-type or $OG10$-type hyperkähler manifolds have been illustrated in \cite{Liu24}. This idea can be used to realize some large-dimensional families given in \cite[Theorem 4.12]{CCC21} that can not be realized by \cite{CKKM19,MW15}. Here we provide an example.

\begin{example}
The moduli space $M_{U(2),\rho_1}$ (for $n=2$) in \cite[Theorem 4.12]{CCC21} is realized in \cite{IKKR17} and a general member $(X,f)$ consists of a moduli space $X$ of stable objects on a triangulated category $\schk$ and a natural involution $f$. According to \cite[Example 3.11]{PPZ23}, the category $\schk$ is equivalent to the derived category $\cate{D}^b(S,\alpha)$ for a K3 surface $(S,H)$ of degree $2$ with a 2-torsion Brauer class $\alpha$ described in \cite[§9.8]{vG05}. Moreover, the moduli space $X$ is isomorphic to $M_{\sigma}(0,H,0)$ for some $\sigma\in\Stab^{\circ}(S,\alpha)$ according to \cite[Theorem 5.1]{CKKM19} and \cite{BM14-MMP} (see also \cite[Proposition 3.3]{MW15} for the untwisted case). The class $\alpha$ is invariant under the natural involution $\pi$ on $S$, so the pushforward $\pi_*$ is an involution on $\cate{D}^b(S,\alpha)$. By description of the transcendental lattice of $(S,\alpha)$ and (the twisted version of) Proposition \ref{K3-invariant-lattices-compatible}, the involution $f$ is induced by $\pi_*$ via sending $[E]$ to $[\pi_*E]$ for any $E\in M_{\sigma}(0,H,0)$.
\end{example}

\section*{Acknowledgment}
This article is inspired by the mini-courses in \emph{SinG-Trento 2024}. The author would like to thank Chiara Camere and Alessandra Sarti, two speakers for the mini-courses, for their interest and comments. Also, the author wants to thank his advisors Laura Pertusi and Paolo Stellari for helpful remarks. In addition, the author wants to thank Daniele Faenzi for exchanges. Moreover, the author wants to thank Arend Bayer, Kuan-Wen Lai, and his grand-advisor Bert van Geemen for answering his questions. Last but not least, the author would like to thank the anonymous referees whose comments improve this article a lot.


\end{document}